
\input amstex
\documentstyle{amsppt}
\NoBlackBoxes
\magnification=1200
\parindent 20 pt
\nologo
\define \bfmtilde{\underset\sim\to {\bold m}}

\define \Atilde{\underset\sim\to A}

\topmatter
\centerline{\bf Killing Luzin and Sierpinski sets} 
\vskip .20in
\centerline{by}
\vskip .20in
\centerline{H. Judah\footnote""{The authors are partially supported by the Basic
Research Foundation, Israel Academy of Science.}} \centerline{Department of
Mathematics and Computer Science} \centerline{Bar-Ilan University} 
\centerline{Ramat-Gan 52900, Israel}
\vskip .20in \centerline{and}\vskip .20in \centerline{S.
Shelah\footnote"*"{publication number 478} } \centerline{Institute of
Mathematics}
\centerline{Hebrew University of Jerusalem}
\centerline{Givat Ram, Jerusalem}
\centerline{Israel}
\centerline{and}
\centerline{Rutgers University}
\vskip .20in
\abstract{We will kill the old Luzin and Sierpinski sets in order to build a
model where $U(\Cal M) = U(\Cal N)=\aleph_1$ and there are neither Luzin nor
Sierpinski sets.
Thus we answer a question of J. Steprans, communicated by S. Todorcevic on
route from Evans to MSRI.} \endabstract
\endtopmatter

\newpage
\baselineskip 20pt 
In this note we will build a model where there are non-measurable sets and
non-meager sets of size $\aleph_1$ and there are neither Luzin nor
Sierpinski sets.
All our notation is standard and can be found in \cite{Ku},
\cite{BJ1}.
Let us start with the basic concept, underlying this work.

Let $U(\Cal M)$ be the minimal cardinal of a non-meager set.

Let $U(\Cal N)$ be the minimal cardinal of a non-null set.

We say that a set of reals $X$  is a Luzin set if $X$ is uncountable and $X\cap
M$ is countable for every meager set $M.$
We say that a set $X$ is a Sierpinski set if $X$ is uncountable and $X\cap N$
is countable, for every null set $N.$

\demo{Fact} (a) If there is a Luzin set, then $U(\Cal M)= \aleph_1.$

(b) If there is a Sierpinski set then $U(\Cal N)=\aleph_1.$\quad $\square$

In \cite{Sh} it was proved that if ZF is consistent, then there is a model where
there are no Luzin sets and $U(\Cal M)=\aleph_1.$
In \cite{BGJS} it was proved that if there is a Sierpinski set, then there is a
non-measurable meager filter on $\omega.$
It was natural to ask if from $U(\Cal N)=\aleph_1$ we can get such a filter.
Clearly it will be enough to answer positively the following question.

(Strepans) Does $U(\Cal N)=\aleph_1$ imply the existency of a Sierpinski set?

We give a negative answer to this question by proving the following
\proclaim{Theorem}
Cons(ZF) $\rightarrow$ Cons(ZFC $+ U(\Cal M)=U(\Cal N) = \aleph_1 +$ there are
neither Luzin nor Sierpinski sets).\endproclaim
We will prove this theorem by iterating with countable support iteration Miller
reals (rational perfect forcing).
We will use the machinery produced by ``preservation theorems'' to show that
the old reals are neither non-meager nor non-measurable sets.
We will show that Miller reals kill Luzin and Sierpinski sets from the ground
model.

The reader can find a complete analysis of Luzin and Sierpsinki sets in
\cite{BJ2}.

\demo{1. Definition}
Let $P = \{T: T\subseteq \omega^{<\omega}\ \&\ T \ \text{is a tree} \ \& \
(\forall s\in T) (s\ \text{is increasing}) \ \& \ (\forall s\in T)(\exists t\in
T)(\exists ^\infty n)(s\subset t^{\wedge}\langle n\rangle\in T)\}.$

Let $\leq$ be defined by $T\leq S$ iff $S\subseteq T.$ \quad $\square$

$\langle P,\leq\rangle$ is called rational perfect forcing (\cite{Mi}) and if
$G\subseteq P$ is generic, then $\bold m=\cap G\in \omega^\omega$ is called a
Miller real.
>From our assumption we have that $\bold m$ is increasing.
\demo{2. Definition}
Let $r\in\omega^\omega$ be increasing.
We define the following set
$$B(r) = \bigcup_{j<\omega}B_j(r)$$
where
$$B_j(r)=\{\eta\in 2^\omega: (\forall i>j)\
(\eta\restriction[r(i),r(i)+10(i+1))\ \text{is not identically zero})\}$$
\demo{3. Fact}
$\mu(B_j(r))\geq 1 - \frac{1}{j+1}$ and $B_j(r)$ is closed. \quad $\square$

\flushpar Therefore $\mu(B(r))=1.$

\proclaim{4. Lemma}
Let $A$ be a set of reals such that $\mu^*(A)>0.$
Let $\bfmtilde$ be the canonical name for the Miller real. Then
$$\Vdash_P``A-B(\bfmtilde)\ \text{is uncountable}".$$
\endproclaim
\demo{Proof}
Let $p\in P,$ \ $p\Vdash_P``\Atilde' = A\backslash B(\bfmtilde) \ \text{is
countable}".$
As $P$ is proper, w.l.o.g.  for some countable set $A^*\subseteq A$ and $q\geq
p$ we have $q\Vdash_P ``\Atilde'\subseteq A^*".$
Let $N\prec\langle H((2^{\aleph_0})^*),\in\rangle$ be countable, $q\in
N,$ \ $A\in N,$ \ $A^*\in N.$
As $\mu^*(A)>0$ there is $\eta\in A,$ \ $\eta$ random over $N.$
Therefore $\eta\notin A^*.$
Let $t\in q$ be such that $mc_q(t)=\{n: t^\wedge\langle n\rangle \in q\}$ is
infinite.
Let us write this set as
$$mc_q(t)=k_\ell^t: \ell<\omega\},$$
where $k^t_\ell<k_{\ell+1}^t.$  
Let $i_t=|t|.$
For $n<\omega,$ we define
$$E_t^n=\{x\in 2^\omega: (\forall \ell\geq n)(x\restriction[k_\ell^t,k_\ell^t
+10(i_t+1))\ \text{is not identically zero})\}.$$
\demo{5. Fact}
$\mu(E_t^n)=0.$ \quad $\square$ 

Therefore $E_t=\bigcup\limits_n E_t^n$ is null and $E_t\in N.$
Therefore $\eta\notin E_t.$
Therefore 
$$D_t=\{k_\ell^t: \eta\restriction[k_\ell^t, k_\ell^t+ 10(i_t+1))\ \text{is
identically zero}\}$$
is infinite.

Now using this we can define, inductively, $q'\geq q$ satisfying
$$\text{if} \ t\in q'\ \text{and} \ mc_q(t)\ \text{is infinite, then}\
mc_{q'}(t)=D_t.$$
Therefore $q'\Vdash_P``\eta\notin B(\bfmtilde)",$ a contradiction. \quad
$\square$ \proclaim{6. Corollary}
If $Y\in V$ is a Sierpinski set, then $Y$ is not a Sierpinski set in
any extension of $V$ containing a Miller real over $V.$ \quad $\square$
\endproclaim
\demo{7. Remark}
The same result can be obtained if you replace Miller real by Laver real.
\demo{8. Definition}
Let $r\in \omega^\omega$ be increasing.
We define the following set
$$T(r) = \bigcup\limits_{j<\omega} [T_j (r)]$$
where $T_j(r)$ is the tree defined by
$$\eta\in [T_j(r)]\ \text{iff}\ \eta\in 2^\omega \ \& \ (\forall i
>j)(\eta(r(i))=0).$$ We say that for a tree $T,$ \ $[T]$ is the set of
$\omega$-branches of $T.$.\quad $\square$
\demo{9. Fact}
$[T_j(r)]$ is a closed nowhere dense set. \quad $\square$

\flushpar Therefore $T(r)$ is a meager set.
\proclaim{10. Lemma}
Let $A$ be a non-meager set of reals.
Let $\bfmtilde$ be the canonical name for the Miller real.
Then
$$\Vdash_P``A\cap T(\bfmtilde) \ \text{is uncountable}".$$
\endproclaim
\demo{Proof}
Let $p\in P$ and let $N\prec(H((2^{\aleph_0})^+),\in)$ be countable such that
$p\in N.$
Then there is $\eta\in A$ such that $\eta$ is Cohen over $N.$
We will find  $q$ such that
$$p\leqq q\in P\ \text{and} \ q\Vdash_P``\eta\in T(\bfmtilde)".$$
Clearly this is enough.
Let $\langle \nu_\rho: \rho\in ^{\omega>}\omega\rangle$ be the list of
splitting nodes of $p$ such that $\rho_1\subsetneq \rho_2$ implies
$\nu_{\rho_{1}}\subsetneq \nu_{\rho_{2}}.$
Thus $\langle \nu_{\rho^\wedge\langle n\rangle}(|\nu_\rho|):
n<\omega\rangle$ are distinct and w.l.o.g. are strictly increasing, so 
$(*)\quad \nu_{\rho^\wedge\langle n\rangle}(|\nu_\rho|)\geq n.$

For each  $\rho\in ^{\omega >}\omega$ let
$$A_\rho=\{n<\omega: \eta\restriction (\text{Range}\ \nu_{\rho^\wedge\langle
n\rangle}\backslash\ \text{Range}\ \nu_\rho)\ \text{is identically zero}\}.$$
\demo{11. Fact}
For $\rho\in ^{\omega>}\omega,$ \ $A_\rho$ is infinite.

\demo{[Proof} $p\in N$ and let $s\in 2^{<\omega}$ be a condition in Cohen
forcing.  
Then there is $n,$  by (*), such that
$$\text{dom}(s)\cap (\text{Range}\ \nu_{\rho^\wedge\langle
n\rangle}\backslash\text{Range} \ \nu_\rho)=\emptyset.$$
Thus we can extend $s$ to $t\in 2^{<\omega}$ such that
$$t\restriction(\text{Range}\ \nu_{\rho^\wedge\langle
n\rangle}\backslash\text{Range}\ \nu_\rho)$$
is identically zero.
Thus, because $\eta$ is Cohen over $N,$ we have that $A_\rho$ is infinite.
\quad $\square$]

\flushpar Now we define $q$ by
$$q = \{\nu\in p: (\forall
\ell\leq|\nu|)(\nu\restriction\ell=\nu_{\rho^\wedge\langle n\rangle} \rightarrow
n\in A_\rho)\}$$ and $q\Vdash_P``\eta\in T(\bfmtilde)".$ \quad
$\square$ \proclaim{12. Corollary}
If $X\in V$ is a Luzin set, then $X$ is not a Luzin set in any extension of $V$
containing a Miller real over $V.$ \quad $\square$ \endproclaim

Now we are ready to show the main Theorem. 
\proclaim{13. Theorem}
$\text{Cons(ZF)}$ implies $\text{Cons(ZFC}+ U(\Cal M)=U(\Cal N)=\aleph_1+$
there are neither Luzin nor Sierpinski sets).
\endproclaim
\demo{Proof}
Let us start with $V=L.$
Let $P_{\omega_{2}}$ be the countable support iteration of $P,$ of
length $\omega_2.$
Then the following holds in $V^{P_{\omega_{2}}}.$
\roster
\item"(i)"
$U(\Cal M)=\aleph_1$: In \cite{Go} it is proved that the property of
being non-meager is preserved by a countable support iteration.
It is easy to see that $P$ satisfies the covering properties established in
\cite{Go}, \S6.20.
Therefore $V\cap 2^\omega$ is a non-meager set in $V^{P_{\omega_{2}}}.$
\item"(ii)"
$U(\Cal N)=\aleph_1$: In \cite{Go} it is proved that the property of
being non-null is preserved by countable support iteration.
In \cite{BJS} it is proved that $P$ satisfies the covering properties
established in \cite{Go} \S6.8.
Therefore $V\cap 2^\omega$ is a non-null set in $V^{P_{\omega_{2}}}.$
\item"(iii)"
There are no Luzin sets in $V^{P_{\omega_{2}}}$: by  corollary 6.
\item"(iv)"
There are no Sierpinski sets in  $V^{P_{\omega_{2}}}$: by  corollary 12. \quad
$\square$
\endroster
\demo{14. Remark}
In the $\omega_2$-iteration of Laver reals we have that $U(\Cal N)=\aleph_1$
and there are no Sierpinski sets. 
We don't know if in this model there are uncountable strongly meager sets.
We know that Miller reals do not kill strong measure zero sets.
This is a consequence of a Rothberger theorem.  See \cite{BJ2}.
\enddemo

\end